\documentclass[a4paper,12pt]{article}
\usepackage{amsmath,amssymb,amsfonts}
\usepackage{hyperref}

\def\mdot{{\cdot}}
\newcommand{\gmod}[1]{#1\hbox{-}\mathsf{Mod}}

\newcommand{\hc}[2]{\Hom_{\mathcal{C}}(#1,#2)}

\newcommand{\ec}[1]{\End_{\mathcal{C}}(#1)}

\def\un{{\bf 1}}
\def\zero{\{0\}}
\def\mpoint{\;\;.}
\def\mvirg{\;\;,}

\def\mpn{\medskip\par\noindent}

\def\normal{\mathop{\underline\triangleleft}}

\def\smp{\smallskip\par}

\def\Id{\hbox{\rm Id}}
\def\Irr{\hbox{\rm Irr}}
\def\Res{\hbox{\rm Res}}
\def\Ind{\hbox{\rm Ind}}
\def\Hom{\hbox{\rm Hom}}
\def\End{\hbox{\rm End}}

\def\Out{\hbox{\rm Out}}
\def\Inf{\hbox{\rm Inf}}
\def\Def{\hbox{\rm Def}}

\def\Iso{\hbox{\rm Iso}}
\def\Im{\hbox{\rm Im}}
\def\Ker{\hbox{\rm Ker}}
\def\Coker{\hbox{\rm Coker}}
\def\Defres{\hbox{\rm Defres}}

\def\Indinf{\hbox{\rm Indinf}}

\def\op{^{op}}

\def\dom{\backslash}
\newcommand{\dirsum}{\mathop{\oplus}}

\newcommand{\flh}[2]{\mathop{\hbox to 4ex{\rightarrowfill}}_{#2}^{#1}\limits}
\newcommand{\flv}[2]{\llap{$#1$}\left\downarrow\vbox to 6mm{}\right.\rlap{$#2$}}

\newcommand{\sur}[1]{\,\overline{\! #1}}
\newcommand{\sumb}[2]{\sum_{{\scriptstyle #1}\atop {\scriptstyle #2}}}

\newcommand{\ressort}[1]{\hskip #1 plus #1 minus #1}

\def\findemo{~\leaders\hbox to 1em{\hss\ \hss}\hfill~\raisebox{.5ex}{\framebox[1ex]{}}\smp}
\newcommand{\carre}[8]{\begin{array}{ccc}
#1&\mathop{\hbox to 12mm{\rightarrowfill}}^{\displaystyle{#2}}\limits&#3\\
\llap{$\displaystyle{#4}$}\left\downarrow\vbox to 6mm{}\right. & & \left\downarrow\vbox to 6mm{}\right.\rlap{$\displaystyle{#5}$}\\
#6&\mathop{\hbox to 12mm{\rightarrowfill}}_{\displaystyle #7}\limits&#8\\
\end{array}}

\newcommand{\carrem}[8]{\begin{array}{ccc}
#1&\mathop{\hbox to 12mm{\rightarrowfill}}^{\displaystyle #2}\limits&#3\\
\llap{$\displaystyle #4$}\left\uparrow\vbox to 6mm{}\right. & & \left\uparrow\vbox to 6mm{}\right.\rlap{$\displaystyle #5$}\\
#6&\mathop{\hbox to 12mm{\rightarrowfill}}_{\displaystyle #7}\limits&#8\\
\end{array}}

\makeatletter
\renewenvironment{enumerate}{\ifnum \@enumdepth >3 \@toodeep\else
      \advance\@enumdepth \@ne
      \edef\@enumctr{enum\romannumeral\the\@enumdepth}\list
      {\csname label\@enumctr\endcsname}{\setlength{\topsep}{1ex}\setlength{\itemsep}{0pt}\usecounter
        {\@enumctr}\def\makelabel##1{\hss\llap{##1}}}\fi}{\endlist}
\renewenvironment{itemize}{\ifnum \@itemdepth >3 \@toodeep\else \advance\@itemdepth \@ne
\edef\@itemitem{labelitem\romannumeral\the\@itemdepth}%
\list{\csname\@itemitem\endcsname}{\setlength{\topsep}{1ex}\setlength{\itemsep}{0pt}\def\makelabel##1{\hss\llap{##1}}}\fi}
{\endlist}

\def\@sect#1#2#3#4#5#6[#7]#8{\ifnum #2>\c@secnumdepth
     \let\@svsec\@empty\else
     \refstepcounter{#1}\edef\@svsec{\csname the#1\endcsname .\hskip .5em}\fi
     \@tempskipa #5\relax
      \ifdim \@tempskipa>\z@
        \begingroup #6\relax
          \@hangfrom{\hskip #3\relax\@svsec}{\interlinepenalty \@M #8\par}%
        \endgroup
       \csname #1mark\endcsname{#7}\addcontentsline
         {toc}{#1}{\ifnum #2>\c@secnumdepth \else
                      \protect\numberline{\csname the#1\endcsname}\fi
                    #7}\else
        \def\@svsechd{#6\hskip #3\relax  
                   \@svsec #8\csname #1mark\endcsname
                      {#7}\addcontentsline
                           {toc}{#1}{\ifnum #2>\c@secnumdepth \else
                             \protect\numberline{\csname the#1\endcsname}\fi
                       #7}}\fi
     \@xsect{#5}}

\def\section{\pagebreak[3]\setcounter{prop}{0}\setcounter{equation}{0}\@startsection{section}{1}{\z@}{6ex plus 6ex}{6ex}{\center\reset@font\large\bf}}
\def\subsection{\pagebreak[3]\refstepcounter{prop}\@startsection{subsection}{2}{\z@}{4ex plus 6ex}{-1em}{\reset@font\bf}}
\def\subsubsection{\@startsection{subsubsection}{3}{\z@}{4ex plus 6ex}{-1em}{\reset@font\it}}
\makeatother

\def\Z{\mathbb{Z}}
\def\1{{1\;\!\!\!{\rm l}}}

\def\Q{\mathbb{Q}}

\def\F{\mathbb{F}}

\def\theprop{\thesection.\arabic{prop}}
\renewenvironment{equation}{\refstepcounter{prop}$$}{\leqno{(\thesection.\arabic{prop}})$$}

\renewenvironment{equation}{\refstepcounter{subsection}\refstepcounter{prop}$$}{\leqno{\bf (\theprop)}$$}

\def\pf{\noindent{\bf Proof: }}
\newenvironment{rem}[1]{\refstepcounter{subsection}\refstepcounter{prop}\mpn{{\bf \thesection.\arabic{prop}.}\ \ \bf#1~:}}{\smp}
\newenvironment{enonce}[1]{\pagebreak[3]\refstepcounter{subsection}\refstepcounter{prop}\mpn{{\bf \thesection.\arabic{prop}.\ \ #1~:}}\begin{it} }{\end{it}\smp}

\def\thesection{\arabic{section}}

\begin{document}
\centerline{\Large\bf Rational $p$-biset functors}
\vspace{.5cm}\par
\centerline{\bf Serge Bouc}
\vspace{1cm}\par
\begin{footnotesize}
{\bf Abstract~:} In this paper, I give several characterizations of {\em rational biset functors over $p$-groups}, which are independent of the knowledge of genetic bases for $p$-groups. I~also introduce a construction of new biset functors from known ones, which is similar to the Yoneda construction for representable functors, and to the Dress construction for Mackey functors, and I show that this construction preserves the class of rational $p$-biset functors.\par
This leads to a characterization of rational $p$-biset functors as additive functors from a specific quotient category of the biset category to abelian groups. Finally, I give a description of the largest rational quotient of the Burnside $p$-biset functor~: when $p$ is odd, this is simply the functor $R_\Q$ of rational representations, but when $p=2$, it is a non split extension of $R_\Q$ by a specific uniserial functor, which happens to be closely related to the functor of units of the Burnside ring.\vspace{.2cm}\par
\end{footnotesize}
\vspace{.5cm}\par
\section{Introduction}\label{intro}
The formalism of biset functors for finite groups yields a unified framework for operations of induction, restriction, inflation, deflation and transport by group isomorphism, which appear in various areas of group representation theory. It was used recently as an essential tool to the solution of some open problems, such as the classification of endo-permutation modules for an arbitrary finite $p$-group (\cite{dadegroup}), or the structure of the group of units of the Burnside ring of such a $p$-group (\cite{burnsideunits}).\par
This formalism seems indeed particularly well suited for studying $p$-groups (where $p$ is a fixed prime number), the main reason for this being that $p$-groups have many normal subgroups, so for these groups, the operations of inflation and deflation play a prominent r\^ole. Also, some important biset functors are only defined on $p$-groups, and cannot be seen in a natural way as the restriction to $p$-groups of a functor defined on arbitrary finite groups~: an example of this situation is the Dade functor for $p$-groups. Such biset functors are called {\em $p$-biset functors}.\par
Among $p$-biset functors, the subclass of {\em rational} functors was introduced in \cite{bisetsections} (see~\ref{defrat}), the name {\em rational} being chosen because these functors behave like the functor of rational representations, in regard to the decomposition of their evaluations at a given $p$-group with respect to {\em faithful elements}. These rational biset functors are generally easy to describe, since they are determined by their knowledge at $p$-groups of normal $p$-rank 1~: this was in particular the key argument to the structure of the group of units of Burnside rings of $p$-groups. Another important property of rational $p$-biset functors is that they form a {\em Serre subclass} of all $p$-biset functors.\par
This paper originates in stimulating and fruitful conversations I had with Laurence Barker and Erg\"un Yal\c c\i n, during my visit at Bilkent University in April 2006. The original question arose when L.~Barker and I understood that we had two a priori different definitions of what a rational $p$-biset functor should be~: L.~Barker (\cite{rhetoric}, see also~\ref{barkerrat}) used a definition akin to the notion of functor on the category $RG$-Morita introduced by I.~Hambleton, L.~R.~Taylor, and E.~B.~Williams in~\cite{htw}. This definition is much more natural that mine, since in particular it applies to arbitrary finite groups, and not only to $p$-groups.\par
So the question was to know whether the two notions could be equivalent. My very first feeling was that the answer should be negative. Also, it seemed clear to  me that rational functors in the sense of L.~Barker should be rational in mine, but the converse was not obvious. Shortly afterwards, independently of each other, and using different methods, E. Yal\c c\i n, L. Barker and I gave a proof of this fact.\par
This paper gives a complete answer to the question of the equivalence of these definitions (Theorem~\ref{le mur}), showing that my first intuition was almost always wrong~: the two definitions of rational $p$-biset functors are equivalent if $p$ is an odd prime number. But for $p=2$, the two definitions are not equivalent, the typical example being the largest rational quotient (in my sense) of the Burnside functor, which is not rational for the second definition. In particular, this definition does not lead to a Serre subclass for~$p=2$.\par

This paper is organized as follows~: in Section~\ref{rappels}, I recall the basic definitions and results on biset functors. In Section~\ref{caractérisation}, I state a necessary and sufficient condition (Theorem~\ref{caract}) for a biset functor to be rational (in my sense), along the lines of Proposition~8.1 of~\cite{dadegroup}. In Section~\ref{Yoneda-Dress}, I introduce a construction for biset functors (\ref{YD}), which is similar both to the Dress construction for Mackey functors, and to Yoneda construction for arbitrary functors. This leads, in Section~\ref{catégorie}, to a characterization (Theorem~\ref{le mur}) of the category of rational biset functors as a subcategory of all biset functors over $p$-groups. In particular, Corollary~\ref{geometric} is a kind of geometric characterization of rational $p$-biset functors. In Section~\ref{conoyau exp}, I describe the largest rational quotient of the Burnside functor~: if $p$ is odd, this functor is equal to $R_\Q$, but if $p=2$, this functor is a non-split extension of $R_\Q$ by a 2-torsion functor which happens to be isomorphic to a subfunctor of the functor of units of the Burnside ring, more precisely the cokernel of the exponential map. At the time, this is just a coincidence, and I could find no natural {\em a~priori} isomorphism between those two objects.
\section{Biset functors}\label{rappels}
\begin{enonce}{Definition}  The {\em biset category}~$\mathcal{C}$ on finite groups is defined as follows~:
\begin{itemize}
\item The objects of $\mathcal{C}$ are finite groups.
\item If $G$ and $H$ are finite groups, then 
$$\Hom_\mathcal{C}(G,H)=B(H\times G\op)$$
is the Burnside group of $(G,H)$-bisets.
\item If $G$, $H$, and $K$ are finite groups, then the composition
$$\Hom_\mathcal{C}(G,H)\times\Hom_\mathcal{C}(H,K)\to\Hom_\mathcal{C}(G,K)$$
is the bilinear extension of the product $(V,U)\mapsto V\times_HU$, where $V$ is a $(K,H)$-biset, and $U$ is an $(H,G)$-biset. Moreover, the identity morphism of the group $G$ in $\mathcal{C}$ is the class of the $(G,G)$-biset $G$, for left and right multiplication.
\end{itemize}
\end{enonce}
A {\em biset functor} is an additive functor from $\mathcal{C}$ to the category $\mathcal{A}b$ of abelian groups. A {\em morphism of biset functors} is a natural transformation of functors. The {\em category of biset functors} is denoted by $\mathcal{F}$.\par
When $p$ is a prime number, the {\em $p$-biset category} $\mathcal{C}_p$ is the full subcategory whose objects are finite $p$-groups. A {\em $p$-biset functor} is an additive functor from $\mathcal{C}_p$ to $\mathcal{A}b$. The {\em category of $p$-biset functors} is denoted by $\mathcal{F}_p$.
\subsection{Elementary morphisms.} This formalism of bisets gives a single framework for the usual operations of induction, restriction, inflation, deflation, and transport by isomorphism via the following correspondences~: 
\begin{itemize} 
\item If $H$ is a subgroup of $G$, then let $\Ind_H^G\in\hc{H}{G}$ denote the set~$G$, with left action of $G$ and right action of $H$ by multiplication.
\item If $H$ is a subgroup of $G$, then let $\Res_H^G\in\hc{G}{H}$ denote the set~$G$, with left action of $H$ and right action of $G$ by multiplication.
\item If $N\normal G$, and $H=G/N$, then let $\Inf_H^G\in\hc{H}{G}$ denote the set~$H$, with left action of $G$ by projection and multiplication, and right action of $H$ by multiplication.
\item If $N\normal G$, and $H=G/N$, then let $\Def_H^G\in\hc{G}{H}$ denote the
 set~$H$, with left action of $H$ by multiplication, and right action of $G$ by  projection and multiplication.
\item If $\varphi: G\to H$ is a group isomorphism, then let $\Iso_G^H=\Iso_G^H(\varphi)\in \hc{G}{H}$ denote the set~$H$, with left action of $H$ by multiplication, and right action of $G$ by taking image by $\varphi$, and then multiplying in~$H$.
\item If $(T,S)$ is a section of the group $G$, i.e. a pair of subgroups of $G$ with $S\normal T\subseteq G$, set
$$\Indinf_{T/S}^G=\Ind_{T}^G\circ\Inf_{T/S}^T\ressort{1cm}\Defres_{T/S}^G=\Def_{T/S}^T\circ\Res_T^G\mpoint$$
Then $\Indinf_{T/S}^G$ is isomorphic to the $(G,T/S)$-biset $G/S$, and $\Defres_{T/S}^G$ is isomorphic to the $(T/S,G)$-biset $S\dom G$.
\end{itemize}
\begin{enonce}{Notation} When $F$ is a biset functor, and $\varphi\in\hc{G}{H}$, the image of the element $u$ of $F(G)$ by the map $F(\varphi)$ will be also be denoted by $\varphi(u)$ or even $\varphi u$.\par
Similarly, the maps $F(\Ind_H^G)$, $F(\Res_H^G)$,\ldots, will be denoted by $\Ind_H^G$, $\Res_H^G$,\ldots, and they will be called induction, restriction,\ldots
\end{enonce}
\subsection{Factorization of morphisms.} Any morphism form $G$ to $H$ in $\mathcal{C}$ is a linear combination with integer coefficients of classes of transitive $(H,G)$-bisets, and isomorphism classes of transitive $(H,G)$-bisets are in one to one correspondence with conjugacy classes of subgroups of the direct product $H\times G$~: the $(H,G)$-biset corresponding to the subgroup $L$ of $H\times G$ is the set of cosets $(H\times G)/L$, for biset action defined by
$$\forall (g,x)\in G^2,\;\forall (h,y)\in H^2,\;\;g.(x,y)L.h=(gx,h^{-1}y)L\mpoint$$
Denote by $p_1(L)$ the projection of $L$ on $H$, and by $p_2(L)$ its projection on $G$. Set moreover
$$k_1(L)=\{h\in H\mid (h,1)\in L\}\ressort{1cm}k_2(L)=\{g\in G\mid (1,g)\in L\}\mpoint$$
Then $k_i(L)\normal p_i(L)$, for $i=1,2$, the group $k_1(L)\times k_2(L)$ is a normal subgroup of $L$, and there are canonical group isomorphisms
$$p_1(L)/k_1(L)\cong L/\big(k_1(L)\times k_2(L)\big)\cong p_2(L)/k_2(L)\mpoint$$
Now (see Lemma~7.4 of~\cite{both}) the biset $(H\times G)/L$ can be factored as
$$\Ind_{p_1(L)}^H\circ\Inf_{p_1(L)/k_1(L)}^{p_1(L)}\circ\Iso(\varphi)\circ\Def_{p_2(L)/k_2(L)}^{p_2(L)}\circ\Res_{p_2(L)}^G\mvirg$$
where $\varphi: p_2(L)/k_2(L)\to p_1(L)/k_1(L)$ is the above isomorphism, sending $gk_2(L)$, for $g\in G$, to $hk_1(L)$, where $h\in H$ is chosen such that $(h,g)\in L$. Note that if $G$ and $H$ are $p$-groups, then this factorization holds within the category $\mathcal{C}_p$, since the groups $p_i(L)$ and $p_i(L)/k_i(L)$ are $p$-groups in this case, for $i=1,2$.
\subsection{Faithful elements.} (see Lemma~2.7 of~\cite{bisetsections})
Let $N$ be a normal subgroup of $G$. Set
$$f_N^G=\sumb{M\normal G}{N\subseteq M}\mu_{\normal G}(N,M)\;\Inf_{G/M}^G\circ\Def_{G/M}^G\mvirg$$
where $\mu_{\normal{G}}$ is the M\"obius function of the poset of normal subgroups of $G$.
This is an element of $\ec{G}$. It is an idempotent, and the elements $f_N^G$, for $N\normal G$, are mutually orthogonal idempotents, whose sum is equal to the identity of $\ec{G}$.\par
If $F$ is a biset functor, set
$$\partial F(G)=F(f_\un^G)F(G)\mpoint$$
It is a direct summand of $F(G)$, called the set of {\em faithful elements} of $F(G)$. It is equal to the intersection of the kernel of all deflation maps $\Def_{G/M}^G$, where~$M$ is a non trivial normal subgroup of $G$.
\begin{rem}{Remark} If $P$ is a $p$-group, then
$$f_\un^P=\sum_{Z\subseteq \Omega_1Z(P)}\mu(\un,Z)\;\Inf_{P/Z}^P\circ\Def_{P/Z}^P\mvirg$$
where $\Omega_1Z(P)$ is the largest elementary abelian subgroup of the center of $P$, and $\mu$ is the M\"obius function of the poset of its subgroups.
\end{rem}
\begin{rem}{Simple functors} Recall some notation and basic results on simple $p$-biset functors : the category of $p$-biset functors is an abelian category. Its simple objects are parametrized by pairs $(H,V)$, where $H$ is a finite $p$-group, and $V$ is a simple $\Z \Out(H)$-module, where $\Out(H)$ is the group of outer automorphisms of $H$ (see e.g.~\cite{doublact}, Proposition~2, page 678, or \cite{both} Proposition~7.10 for details). All the simple functors appearing in the present paper (in Section~\ref{conoyau exp}) correspond to the case where~$V$ is the quotient $\Z/2\Z$ (with obviously trivial $\Out(H)$-action). The corresponding simple functor will be denoted by~$S_{H,\F_2}$.
\end{rem}

\begin{enonce}{Genetic subgroups of $p$-groups} Let $P$ be a $p$-group. A subgroup $Q$ of $P$ is called {\em genetic} if the two following conditions are fulfilled~:
\begin{itemize}
\item The group $N_P(Q)/Q$ has {\em normal $p$-rank 1}, i.e. it is cyclic, generalized quaternion, or dihedral or semi-dihedral of order at least 16.
\item If $x\in P$, then $Q^x\cap Z_P(Q)\subseteq Q$ if and only if $Q^x=Q$, where $Z_P(Q)$ is the subgroup of $N_P(Q)$ defined by $Z_P(Q)/Q=Z\big(N_P(Q)/Q\big)$.
\end{itemize}
\end{enonce}
Two genetic subgroups $Q$ and $R$ of $P$ are said to be {\em linked modulo $P$} if there exist $x$ and $y$ in $P$ such that
$$Q^x\cap Z_P(R)\subseteq R\ressort{1cm}{\rm and}\ressort{1cm}R^y\cap Z_P(Q)\subseteq Q\mpoint$$
This is an equivalence relation on the set of genetic subgroups of $P$. A {\em genetic basis} of $P$ is a set of representatives of equivalence classes for this relation.
\subsection{Rational $p$-biset functors.} \label{defrat}(see Section~7 of \cite{bisetsections})  Let $P$ be a $p$ group, and $\mathcal{G}$ be a genetic basis of $P$. It was shown in Theorem~3.2. of~\cite{bisetsections} that the map
$$\mathcal{I}_\mathcal{G}:\dirsum_{Q\in\mathcal{G}}\Indinf_{N_P(Q)/Q}^P:\dirsum_{Q\in\mathcal{G}}\partial F\big(N_P(Q)/Q\big)\to F(P)$$
is split injective (with explicit left inverse). The functor $F$ is called {\em rational} if for any $p$-group $P$ and any genetic basis $\mathcal{G}$ of $P$, the map $\mathcal{I}_\mathcal{G}$ is an isomorphism. Equivalently, for any $p$-group $P$, there exists a genetic basis with this property.
\section{A characterization of rational $p$-biset functors}\label{caractérisation}
\begin{enonce}{Theorem} \label{caract}Let $F$ be a $p$-biset functor. Then $F$ is rational if and only if the following two conditions are satisfied~:
\begin{enumerate}
\item[(i)] If $P$ is a $p$-group with non cyclic center, then $\partial F(P)=\zero$.
\item[(ii)] If $P$ is any $p$-group, if $E\normal P$ is a normal elementary abelian subgroup of rank~2, if $Z$ is a central subgroup of order $p$ of $P$ contained in $E$, then the map
$$\Res_{C_P(E)}^P\dirsum \Def_{P/Z}^P : F(P)\longrightarrow F\big(C_P(E)\big)\dirsum F(P/Z)$$
is injective.
\end{enumerate}
\end{enonce}
\pf Suppose first that $F$ is rational. Let $P$ be a $p$-group, and $\mathcal{G}$ be a genetic basis of $P$. Since $F$ is rational, the map
$$\mathcal{D}_\mathcal{G}: F(P)\flh{\dirsum b_Q^P}{}\dirsum_{Q\in\mathcal{G}}\partial F\big(N_P(Q)/Q\big)$$
is injective. Recall that the component $b_Q^P$ of the map $\mathcal{D}_\mathcal{G}$ indexed by the genetic subgroup $Q\in\mathcal{G}-\{P\}$ is equal to $\Defres_{N_P(Q)/Q}^P-\Defres_{N_P(Q)/\hat{Q}}^P$, where $\hat{Q}$ is the subgroup defined by $\hat{Q}/G=\Omega_1Z\big(N_P(Q)/Q\big)$. If $Q=P$, then the corresponding component $b_P^P$ is equal to $\Def_{P/P}^P$.\par
Suppose that the center $Z(P)$ of $P$ is not cyclic. Then for any $Q\in\mathcal G$, the group $Z=Q\cap Z(P)$ is non trivial~: indeed, the group $Z(P)/Z$ is isomorphic to a subgroup of the center of $N_P(Q)/Q$, which is cyclic since $N_P(Q)/Q$ has normal $p$-rank~1. Now for $Q\neq P$, the maps $\Defres_{N_P(Q)/Q}^P$ and $\Defres_{N_P(Q)/\hat{Q}}^P$ factor as
\begin{eqnarray*}
\Defres_{N_P(Q)/Q}^P&=&\Defres_{N_{\sur{P}}(\sur{Q})/\sur{Q}}^{\sur{P}}\Def_{P/Z}^P\\
\Defres_{N_P(Q)/\hat{Q}}^P&=&\Defres_{N_{\sur{P}}(\sur{Q})/\sur{\hat{Q}}}^{\sur{P}}\Def_{P/Z}^P\\
\end{eqnarray*}
where $\sur{P}=P/Z$, where $\sur{Q}=Q/Z$, and $\sur{\hat{Q}}=\hat{Q}/Z$. Moreover if $Q=P$, then $\Def_{P/P}^P=\Def_{\sur{P}/\sur{P}}^{\sur{P}}\Def_{P/Z}^P$. It follows that if $u\in \partial F(P)$, then $\mathcal{D}_\mathcal{G}(u)=0$, since $u$ is mapped to 0 by the proper deflation $\Def_{P/Z}^P$. Thus $u=0$, and $\partial F(P)=\zero$ when the center of $P$ is non cyclic. So if $F$ is rational, then it satisfies condition~$(i)$ of the theorem.\par
For condition $(ii)$, there is nothing to check if $E$ is central in $P$, since in that case the map $\Res_{C_P(E)}^P$ is the identity map, hence it is injective. And if~$E$ is not central in $P$, then $Z=E\cap Z(P)$. In this situation, by Corollary~5.3 of~\cite{dadegroup}, there exists a genetic basis $\mathcal{G}$ which admits a decomposition $\mathcal{G}=\mathcal{G}_1\sqcup\mathcal{G}_2$ as a disjoint union, where 
\begin{equation}\label{condbase}
\begin{minipage}{0.901\textwidth}\begin{itemize}
\item $\mathcal{G}_2=\{Q\in\mathcal{G}\mid Q\supseteq Z\}$, and
\item the set $\mathcal{H}_1=\{{^xQ}\mid Q\in\mathcal{G}_1,\;x\in [P/C_P(E)]\}$ is the subset of a genetic basis of $C_P(E)$ consisting of subgroups not containing~$Z$.
\end{itemize}
\end{minipage}
\end{equation}
Let $u\in F(P)$ be an element such that $\Res_{C_P(E)}^Pu=0$ and $\Def_{P/Z}^Pu=0$. If $Q\in\mathcal{G}_2$, then $b_Q^{P}=b_{\sur{Q}}^{\sur{P}}\circ \Def_{P/Z}^P$, where $\sur{P}=P/Z$ and $\sur{Q}=Q/Z$. It follows that $b_Q^P(u)=0$. Similarly if $Q\in \mathcal{G}_1$, then $N_P(Q)\subseteq C_P(E)$, and $b_Q^P=b_Q^{C_P(E)}\circ \Res_{C_P(E)}^P$, so $b_Q^P(u)=0$ again. Hence $\mathcal{D}_\mathcal{G}(u)=0$, and $u=0$ since $\mathcal{D}_\mathcal{G}$ is injective when~$F$ is rational. Thus $F$ satisfies $(ii)$.\par
Conversely, suppose that $F$ is a $p$-biset functor satisfying $(i)$ and $(ii)$. Showing that $F$ is rational amounts to showing that the map $\mathcal{D}_\mathcal{G}$ is injective, for any $p$-group $P$ and any genetic basis $\mathcal{G}$ of $P$. Actually, by Lemma~7.3 of~\cite{bisetsections}, it is enough to show that $\mathcal{D}_\mathcal{G}$ is injective, for any $p$-group $P$ and for one particular genetic basis $\mathcal{G}$ of $P$. I will prove this by induction on the order of $P$.\par
If $P$ is the trivial group, then the map $\mathcal{D}_\mathcal{G}$ is the identity map, so there is nothing to prove in this case, and this starts induction. Suppose that $P$ is a $p$-group, such that the map $\mathcal{D}_\mathcal{H}$ is injective for any $p$-group $P'$ with $|P'|<|P|$ and any genetic basis $\mathcal{H}$ of $P'$. Let $\mathcal{G}$ be a genetic basis of $P$, and $u\in \Ker\,\mathcal{D}_\mathcal{G}$. Let $N$ be a non-trivial normal subgroup of $P$. Then the set $\mathcal{H}=\{Q/N\mid Q\in\mathcal{G},\;\supseteq N\}$ is a genetic basis of $P/N$. Moreover one checks easily that $b_{\sur{Q}}^{\sur{P}}\circ\Def_{P/N}^P=b_{Q}^P$ for $Q\supseteq N$, where overlines denote quotients by~$N$, as above. Hence $\mathcal{D}_\mathcal{H}(\Def_{P/N}^Pu)=0$, and by induction hypothesis, it follows that $\Def_{P/N}^Pu=0$, and $u\in\partial F(P)$. In other words $\Ker\,\mathcal{D}_\mathcal{G}\subseteq \partial F(P)$.\par
Hence if the center of $P$ is not cyclic, then condition $(i)$ implies that $\Ker\,\mathcal{D}_\mathcal{G}=\zero$, and I can suppose that the center of $P$ is cyclic.\par
If $P$ has normal $p$-rank~1, then the trivial subgroup is a genetic normal subgroup of $P$, hence is belongs to every genetic basis of~$P$. The corresponding map $b_\un^P: F(P)\to \partial F(P)$ is equal to $F(f_\un^P)$. Hence it is equal to the identity map on $\partial F(P)$. Thus $u=F(f_\un^P)(u)=b_\un^P(u)=0$, and the map $\mathcal{D}_\mathcal{G}$ is injective in this case.\par
Finally I can suppose that $P$ admits a normal and non central subgroup~$E$ which is elementary abelian of rank~2. In this case the group $Z=E\cap Z(P)$ has order $p$, and by Corollary~5.3 of~\cite{dadegroup}, there exists a genetic basis $\mathcal{G}$ of $P$ with the two properties~(\ref{condbase}) above. With the same notation, let $\mathcal{H}$ be a genetic basis of $C_P(E)$ containing $\mathcal{H}_1$. Set $v=\Res_{C_P(E)}^Pu$. Then for any $Q\in \mathcal{G}_1$ and any $x\in [P/C_P(E)]$, one has $N_P(Q)\subseteq C_P(E)$ by Lemma~5.2 of~\cite{dadegroup}, and it follows that
\begin{eqnarray*}
b_{^xQ}^{C_P(E)}(v)&=&b_{^xQ}^{C_P(E)}\Res_{C_P(E)}^Pu\\
&=&\Defres_{N_P({^xQ})/{^xQ}}^Pu-\Defres_{N_P({^xQ})/{^x\hat{Q}}}^Pu\\
&=&^x\big(\Defres_{N_P({Q})/{Q}}^Pu-\Defres_{N_P({Q})/{\hat{Q}}}^Pu\big)\\
&=&^xb_Q^Pu=0\;\;,
\end{eqnarray*}
since $^xu=u$ for $x\in P$ and $u\in F(P)$. Thus $b_R^{C_P(E)}(v)=0$ for any $R\in \mathcal{H}_1$. And if $R\in \mathcal{H}-\mathcal{H}_1$, then $R\supseteq Z$, and
$$b_R^{C_P(E)}(v)=b_{R/Z}^{C_P(E)/Z}\Res_{C_P(E)/Z}^{P/Z}\Def_{P/Z}^Pu=0\;\;.$$
It follows that $\mathcal{D}_\mathcal{H}(v)=0$, hence $v=0$  by induction hypothesis. Now $\Res_{C_P(E)}^Pu=0$ and $\Def_{P/Z}^Pu=0$. Condition $(ii)$ implies that $u=0$, so the map $\mathcal{D}_\mathcal{G}$ is injective.\findemo
\begin{rem}{Remark} The argument in the last part of the proof is essentially the same as the one used in the proof of Theorem~8.2 of~\cite{dadegroup}.
\end{rem}
Recall some notation from~\cite{dadegroup}~:
\begin{enonce}{Notation} \begin{itemize}\label{rappeldade}
\item Let $E$ be an elementary abelian $p$-group of rank 2. Define an element~$\varepsilon_E$ of $B(E)$ by
$$\varepsilon=\varepsilon_E=E/\un-\sumb{ F\subseteq E}{|F|=p}E/F+pE/E\mpoint$$
Denote by $B_\varepsilon$ the subfunctor of $B$ generated by $\varepsilon$, i.e. the smallest subfunctor~$F$ of $B$ such that $F(E)\ni \varepsilon$. Equivalently, for any $p$-group~$P$
$$B_\varepsilon(P)=\Hom_{\mathcal{C}_p}(E,P)(\varepsilon)\mpoint$$

\item Let $X$ be either an extraspecial group of order $p^3$ and exponent~$p$ if~$p$ is odd, or a dihedral group of order 8 if~$p=2$. Let $Z$ denote the center of $X$. Let $I$ and $J$ be two subgroups of order $p$ in $X$, different from $Z$, and not conjugate in $X$. Define an element $\delta$ of $B(X)$ by
$$\delta=\delta_{X,I,J}=(X/I-X/IZ)-(X/J-X/JZ)\mpoint$$
Denote by $B_\delta$ the subfunctor of $B$ generated by $\delta$. Thus, for any $p$-group~$P$
$$B_\delta(P)=\Hom_{\mathcal{C}_p}(X,P)(\delta)\mpoint$$
\item Let $K$ denote the kernel of the natural transformation $B\to R_\Q$. Thus~$K$ is a biset subfunctor of $B$.
\end{itemize}
\end{enonce}
\begin{rem}{Remark} These elements $\varepsilon$ and $\delta$ have been introduced in Corollary~6.5 and Notation~6.9 of~\cite{dadegroup}. In the case $p=2$, the element $\delta$ was denoted by $\delta_3$ there. It is easy to check that $\varepsilon_E\in K(E)$ and that $\delta\in K(X)$. So $B_\varepsilon$ and $B_\delta$ are also the subfunctors of $K$ generated by $\varepsilon$ and $\delta$ respectively, that is $B_\varepsilon=K_\varepsilon$ and $B_\delta=K_\delta$. Moreover Lemma~6.11 of~\cite{dadegroup} shows that $\Res_{JZ}^X\delta=\varepsilon_{JZ}$, so in particular $B_\varepsilon\subseteq B_\delta$.
\end{rem}
\begin{rem}{Remark} If $E$ and $E'$ are elementary abelian $p$-groups of rank 2, and if~$\varphi$ is any group isomorphism from $E$ to $E'$, then $\varphi(\varepsilon_E)=\varepsilon_{E'}$. This shows that the subfunctor $B_\varepsilon$ does not depend on the choice of a particular group~$E$. Similarly, if~$X'$ is another extraspecial group of order $p^3$ and exponent $p$, when~$p$ is odd, or another dihedral group of order 8, when~$p=2$, and if~$I'$ and $J'$ are two non conjugate non central subgroups of order $p$ of $X'$, then there exists a group isomorphism $\varphi~:X\to X'$ such that $\varphi(I)=I'$ and $\varphi(J)=J'$. Thus $\varphi(\delta_{X,I,J})=\delta_{X',I',J'}$, and it follows that $B_\delta$ does not depend on the choice of the groups $X$, $I$ and $J$.
\end{rem}
It was shown in Theorem~6.12 of~\cite{dadegroup} that if~$p$ is odd, then $K=B_\delta$. The following provides an alternative proof for this theorem~:
\begin{enonce}{Theorem}\label{brat} The functor $B/B_\delta$ is rational.
\end{enonce}
\pf Let $A$ be any injective cogenerator for the category of $\Z$-modules (e.g. $A=\Q/\Z$). Recall that the $A$-dual $F^0$ of a biset functor $F$ is defined by
$$F^0(P)=\Hom_\Z(F(P),A)\ressort{1cm}F^0(\varphi)={^tF(\varphi\op)}$$
for any $p$-group $P$ and any morphism $\varphi$ in the category $\mathcal{C}_p$. Since $A$ is an injective cogenerator, the canonical map $F\to (F^0)^0$ is injective. Taking $F=B/B_\delta$, I will show that $F^0$ is rational. By Proposition~7.4 of~\cite{bisetsections}, it will follow that $(F^0)^0$ and any subfunctor of it are rational, so $F$ will be rational.\par
Suppose first that $P$ is a $p$-group with non cyclic center. Then Lemma~6.8 of~\cite{dadegroup} shows that $\partial B(P)\subseteq B_\varepsilon(P)\subseteq B_\delta(P)$. Hence $\partial B/B_\delta(P)=\zero$. Since $\partial F^0(P)\cong\Hom_\Z(\partial F(P),A)$, it follows that $\partial F^0(P)=\zero$, so condition~$(i)$ of Theorem~\ref{caract} is satisfied by $F^0$.\par
Now let $P$ be any $p$-group, and let $E$ be a normal subgroup of $P$, which is elementary abelian of rank~2. Let moreover $Z$ be a central subgroup of order~$p$ of~$P$, contained in $E$. Let $f\in F^0(P)$, such that $\Res_{C_P(E)}^P(f)=0$ and $\Def_{P/Z}^Pf=0$. Checking condition $(ii)$ of Theorem~\ref{caract} amounts to showing that $f=0$.\par
If $E$ is central in $P$, there is nothing to do, since $C_P(E)=P$ and $\Res_{C_P(E)}^P(f)=f$. So I can assume that $C_P(E)$ has index $p$ in $P$. Let $\tilde{f}~:B(P)\to I$ be the homomorphism obtained by composing $f$ with the projection $B(P)\to B/B_\delta(P)$. The hypothesis on $f$ means that $\tilde{f}(P/Q)=0$ if $Q\subseteq C_P(E)$ or if $Q\supseteq Z$. Let $Q$ be any subgroup of $P$ with $Q\not\subseteq C_P(E)$ and $Q\not\supseteq Z$. Then $Q\not\supseteq E$. Moreover $QC_P(E)=P$. If the intersection $Q\cap E$ has order~$p$, it is normalized, hence centralized by $Q$, and also by $C_P(E)$, hence $Q\cap E\subseteq Z(P)$. Since $Q\cap E\neq Z$, it follows that $E=Z\mdot(Q\cap E)$ is central in $P$, a contradiction. Hence $Q\cap E=\un$.\par
Now the group $C_Q(E)=Q\cap C_P(E)$ is a normal subgroup of $Q E$. The factor group $R=QE/C_Q(E)$ has order $p^3$, and is the semidirect product of the group $\sur{Q}=Q/C_Q(E)$, of order $p$, by the normal elementary abelian subgroup group $\sur{E}=EC_Q(E)/C_Q(E)\cong E$. Moreover $R$ is not abelian, since otherwise $[Q,E]\subseteq C_Q(E)\cap E=\un$, so $Q\subseteq C_P(E)$. It follows that $R$ is extraspecial of order $p^3$ and exponent $p$ if~$p$ is odd, or dihedral of order 8 if~$p=2$.\par
The center of $R$ is the group $\sur{Z}=ZC_Q(E)/C_Q(E)$. The group $\sur{Q}$ is a non central subgroup of order $p$ of $R$. If $F$ is any subgroup of order $p$ of $E$, different from $Z$, then $\sur{F}=FC_Q(E)/C_Q(E)$ is another non central subgroup of order $p$ of $R$, not conjugate to $\sur{Q}$ in $R$. So in the group $R$, the element 
$$d=(R/\sur{Q}-R/\sur{Q}\sur{Z})-(R/\sur{F}-R/\sur{E})$$
is in $B_\delta(R)$, since it is the image of $\delta$ by a suitable group isomorphism. Taking induction-inflation from the section $(QE,C_Q(E))$ of $P$ gives the element
$$\tilde{d}=\Indinf_{QE/C_Q(E)}^Pd=(P/Q-P/QZ)-\big(P/FC_Q(E)-P/EC_Q(E)\big)\mvirg$$
of $B_\delta(P)$. Now $\tilde{f}$ vanishes on $B_\delta(P)$, thus $\tilde{f}(\tilde{d})=0$, and
$$\tilde{f}(P/Q)=\tilde{f}(P/QZ)+\tilde{f}\big(P/FC_Q(E)\big)-\tilde{f}\big(P/EC_Q(E)\big)=0\mvirg$$
since $QZ\supseteq Z$ and since $FC_Q(E)\subseteq EC_Q(E)\subseteq C_P(E)$. \par
This shows that $\tilde{f}=0$, hence $f=0$, and the functor $F^0$ satisfies both conditions of Theorem~\ref{caract}. Hence $F^0$ is rational, and $F$ is rational as well, since it is a subfunctor of the $A$-dual of $F^0$.\findemo
\begin{enonce}{Corollary}{\rm [\cite{dadegroup}, Theorem~6.12-1]}\label{th sur K} If $p$ is odd, then $K=B_\delta$.
\end{enonce}
\pf Indeed, since $B/B_\delta$ is rational, its subfunctor $K/B_\delta$ is also rational. Then for any $p$-group $P$, and for any genetic basis $\mathcal{G}$ of $P$, the map
$$\dirsum_{Q\in\mathcal{G}}\limits\Indinf_{N_P(Q)/Q}^P~:~\dirsum_{Q\in\mathcal{G}}\limits\partial (K/B_\delta)\big(N_P(Q)/Q\big)\to (K/B_\delta)(P)$$
is an isomorphism. But the groups $N_P(Q)/Q$ have normal $p$-rank~1, hence they are cyclic since~$p$ is odd. Now if~$R$ is cyclic, then $K(R)=\zero$. It follows that $K/B_\delta(P)=\zero$, hence $K=B_\delta$.\findemo
\begin{rem}{Remark} The structure of the functor $K/B_\delta$ in the case $p=2$ will be described in Section~\ref{conoyau exp}.
\end{rem}
\begin{rem}{Remark} The key point in the proof of Theorem~\ref{brat} is the construction of a section $(T,S)$ of $P$ such that $T/S$ is extraspecial of order~$p^3$ and exponent~$p$, or dihedral of order~8. The same idea was used in Lemma~6.15 of~\cite{dadegroup}.
\end{rem}
\section{Yoneda-Dress constructions and rational functors}\label{Yoneda-Dress}
In this section I will introduce a construction which is connected both to the Yoneda functors in the case of the biset category $\mathcal{C}_p$, and to the Dress construction for Mackey functors (see \cite{dress}, or~\cite{green}~1.2).
\begin{enonce}{Notation} Let $H$ be a finite $p$-group. Let $\pi_H~:\mathcal{C}_p\to \mathcal{C}_p$  denote the correspondence sending a $p$-group $P$ to the direct product $P\times H$. If $P$ and $Q$ are $p$-groups, and if $U$ is a $(Q,P)$-biset, denote by $U\times H$ the $(Q\times H,P\times H)$-biset equal to the cartesian product of $U$ and $H$ as a set, with double action defined by
$$\forall y\in Q,\forall x\in P,\forall (h,k,l)\in H^3,\;(y,l)(u,h)(x,k)=(yux,lhk)\mpoint$$
This correspondence extends linearly to an additive map 
$$B(Q\times P\op)\to B\big((Q\times H)\times (P\times H)\op\big)$$
also denoted by $f\mapsto\pi_H(f)$.
\end{enonce}
\begin{enonce}{Lemma} The correspondence $\pi_H$ is an endofunctor of~$\mathcal{C}_p$, which is additive on morphisms.
\end{enonce}
\pf This is straightforward.\findemo
\begin{enonce}{Notation} \label{YD}If $H$ is a $p$-group, and $F$ is a $p$-biset functor, I denote by~$F_H$ the $p$-biset functor obtained by precomposition with $\pi_H$.
\end{enonce}
In other words, for any $p$-group $P$, one has $F_H(P)=F(P\times H)$. If $U$ is a $(Q,P)$-biset, then $\pi_H(U)=F(U\times H)$.
\begin{enonce}{Lemma}\label{f_1^P étendu} Let $P$, $Q$ and $H$ be $p$-groups. If $U$ is a $(P,P)$-biset and $V$ is a $(P\times H,Q)$ biset, then there is an isomorphism of $(P\times H,Q)$-bisets
$$(U\times H)\times_{(P\times H)}V\cong U\times_PV\mvirg$$
where the double action on $U\times_PV$ is induced by
$$(x,h)(u,v)q=(xu,hvq)\mvirg$$
for $x\in P$, $h\in H$, $y\in Q$, $u\in U$, and $v\in V$.\par
In particular, if for any $v\in V$, the intersection of the stabilizer of $v$ in~$P\times H$ with $Z(P)\times\{1\}$ is non trivial, then 
$$(f_\un^P\times H)\times_{(P\times H)}V=0\mpoint$$
\end{enonce}
\pf The first assertion is straightforward. The second one follows from the fact that $V$ is a union of $P$-sets which are inflated from proper quotients of $P$.\findemo
\begin{enonce}{Proposition} \label{Yoneda-Dress rationnel}Let $H$ be a $p$-group, and let $F$ be a rational $p$-biset functor. Then $F_H$ is rational.
\end{enonce}
\pf To prove the proposition, I will check that $F_H$ fulfills the conditions $(i)$ and $(ii)$ of Theorem~\ref{caract}.\par
Let $P$ be a $p$-group with non cyclic center. Let $Q$ be a genetic subgroup of $P\times H$. Since $Z(P)\times\{1\}\subseteq Z(P\times H)$, the group $Z(P)\times\{1\}$ is contained in the center of $N_{P\times H}(Q)$, hence it maps to the center of $N_{P\times H}(Q)/Q$, which is cyclic. Thus if the center of $P$ is not cyclic, then $Z(P)\times\{1\}$ intersects $Q$ non trivially. In other words $Z(P)\cap k_1(Q)\neq \un$.\par
Now, setting $R=P\times H$, there are two cases~: either $Q=R$. In this case~$\gamma_Q$ is just and $(R,R)$ biset of cardinality~1. Thus
$$(f_\un^P\times H)\times_R\gamma_Q=f_\un^P\times_P(P/P)=0\mvirg$$
since $P$ is non trivial. And if $Q\neq R$, then
$$\gamma_Q=(R\times R)/\Delta_{N_R(Q),Q}-(R\times R)/\Delta_{N_R(Q),\hat{Q}}\mvirg$$
where $\hat{Q}$ is defined by the condition that $\hat{Q}/Q$ is the only subgroup of order~$p$ in the center of $N_R(Q)/Q$. The left stabilizer in $R$ of the point $\Delta_{N_R(Q),Q}$ of $(R\times R)/\Delta_{N_R(Q),Q}$ is equal to $Q$, which intersects $Z(P)\times\{1\}$ non trivially by the previous remarks. Similarly, the  left stabilizer in $R$ of the point $\Delta_{N_R(Q),\hat{Q}}$ of $(R\times R)/\Delta_{N_R(Q),\hat{Q}}$ is equal to $\hat{Q}$, which also intersects $Z(P)\times\{1\}$ non trivially. By Lemma~\ref{f_1^P étendu}, it follows that 
$$(f_\un^P\times H)\times_{(P\times H)}\gamma_Q=0\mvirg$$
for any genetic subgroup $Q$ of $R$. But since $F$ is rational, the sum of idempotents $\gamma_Q$, for $Q$ in a genetic basis of $P\times H$, is the identity map of $F(P\times H)$. Thus $f_\un^P\times H$ acts by 0 on $F(P\times H)$. In other words, the idempotent $f_\un^P$ acts by 0 on $F_H(P)$. Hence $\partial F_H(P)=0$, and $F_H$ fulfills condition $(i)$ of Theorem~\ref{caract}.\par
For condition $(ii)$, suppose that $E$ is a normal subgroup of $P$, and that~$E$ is elementary abelian of rank~2. Let $Z$ be a central subgroup of order~$p$ in~$P$, contained in~$E$. It is easy to see that the restriction map from $F_H(P)$ to $F_H\big(C_P(E)\big)$ is equal to the restriction map $\Res_{C_P(E)\times H}^{P\times H}$, after the identifications $F_H(P)=F(P\times H)$ and $F_H\big(C_P(E)\big)=F\big(C_P(E)\times H\big)$. Similarly, the deflation map from $F_H(P)$ to $F_H(P/Z)$ is equal to the map $\Def_{(P\times H)/(Z\times\{1\})}^{P\times H}$ from $F(P\times H)$ to $F\big((P/Z)\times H\big)$, after the identification $(P\times H)/(Z\times\{1\})=(P/Z)\times H$. Now $E'=E\times\{1\}$ is a normal subgroup of $P\times H$, which is elementary abelian of rank 2. Its centralizer in $P\times H$ is $C_P(E)\times H$. Moreover $E'$ contains the subgroup $Z'=Z\times\{1\}$, which is central of order $p$ in $P\times H$. Condition $(ii)$ for the functor $F$ implies that
$$\Ker\,\Res_{C_P(E)\times H}^{P\times H}\cap \Ker\,\Def_{(P/Z)\times H}^{P\times H}=\zero\mpoint$$
Thus $F_H$ fulfills condition $(ii)$ of Theorem~\ref{caract}, and $F_H$ is rational.\findemo
\section{The category of rational $p$-biset functors}\label{catégorie}
In this section, I will show that the rational $p$-biset functors are exactly the additive functors $\mathcal{C}_{p}\to \gmod{\Z}$ which factor through the quotient category $\mathcal{C}_{p}/\delta$ obtained by factoring out every morphism $f: P\to Q$ in $\mathcal{C}_{p}$ which lies in $B_\delta(Q\times P\op)$. The first thing to check for doing this is that this class of morphisms is a two sided ideal in $\mathcal{C}_{p}$~:
\begin{enonce}{Lemma} Let $P$, $Q$, $R$, $S$ be $p$-groups. Then
$$B(S\times Q\op)\times_QB_\delta(Q\times P\op)\times_PB(P\times R\op)\subseteq B_\delta(S\times R\op)\mpoint$$
\end{enonce}
\pf Since $B_\delta(Q\times P\op)=\Hom_{\mathcal{C}_p}(P,Q)(\delta)$, an element in $B_\delta(Q\times P\op)$ is linear combination of isomorphism classes of $(Q,P)$-bisets of the form $U\times_X\delta$, where $U$ is some $(Q\times P\op,X)$-biset, i.e. a $(Q,P\times X)$-biset. Now if $V$ is an $(S,Q)$-biset, and $W$ is a $(P,R)$ biset, then
$$V\times_Q(U\times_X\delta)\times_PW\cong (V\times_QU\times_PW)\times_X\delta$$
in $B(S\times R\op)$, where the $(S, R\times X)$-biset structure on $V\times_QU\times_PW$ is induced by
$$s(v,u,w)(r,x)=(sv,u(1,x),wr)\mvirg$$
for $s\in S$, $r\in R$, $v\in V$, $u\in U$, $w\in W$, and $x\in X$. The lemma follows.\findemo
\begin{enonce}{Definition} Let $\mathcal{C}_{p}/\delta$ denote the following category~: \begin{itemize}
\item The objects are finite $p$-groups.
\item If $P$ and $Q$ are finite $p$-groups, then 
$$\Hom_{\mathcal{C}_{p}/\delta}(P,Q)=B(Q\times P\op)/B_\delta(Q\times P\op)\mpoint$$
\item Composition in $\mathcal{C}_{p}/\delta$ is induced by the product
$$-\times_Q-~: B(R\times Q\op)\times B(Q\times P\op)\to B(R\times P\op)$$
defining the composition in $\mathcal{C}_{p}$.
\end{itemize}
Let $\rho_\delta$ denote the natural functor from $\mathcal{C}_{p}$ to $\mathcal{C}_{p}/\delta$~: so $\rho_\delta$ is equal to the identity on objects, and to the natural projection on morphisms.
\end{enonce}
\begin{enonce}{Theorem} \label{le mur}Let $F$ be a $p$-biset functor. The following conditions are equivalent~:
\begin{enumerate}
\item[(i)] The functor $F$ is rational.
\item[(ii)] The functor $F$ factors through $\rho_\delta$. In other words, if $f\in \Hom_{\mathcal{C}_{p}}(P,Q)$ lies in $B_\delta(Q\times P\op)$, then $F(f)=0$.
\end{enumerate}
\end{enonce}
\pf Suppose first that $(ii)$ holds. Let $\mathcal{S}$ be a set of representatives of isomorphism classes of finite $p$-groups. For any $P\in\mathcal{S}$, choose a set $\Gamma_P$ of generators of $F(P)$ as an abelian group. By Yoneda lemma, the set of natural transformations from $B_P$ to $F$ is in one to one correspondence with $F(P)$, and this correspondence is as follows~: if $s\in F(P)$, then for any $p$-group $Q$, define a map 
$$\sigma_{P,s,Q}~:B_P(Q)=B(Q\times P\op)\to F(Q)$$
by $\sigma_{P,s,Q}(f)=f(s)$. Then the maps $\sigma_{P,s,Q}$ define a morphism of functors $\sigma_{P,s}~: B_P\to F$. If condition $(ii)$ holds, this morphism of functors factors through the quotient functor $(B/B_\delta)_P$, which is rational by Theorem~\ref{brat} and Proposition~\ref{Yoneda-Dress rationnel}. This gives a morphism
$$\dirsum_{P\in\mathcal{S}}\dirsum_{s\in\Gamma_P}(B/B_\delta)_P\to F\mvirg$$
which is obviously surjective. Now any direct sum of rational $p$-biset functors is a rational $p$-biset functor, and any quotient of a rational $p$-biset functor is rational. So $F$ is rational, and $(i)$ holds.

Conversely, suppose now that $F$ is rational. Let $P$ and $Q$ be $p$-groups. Any morphism in $B_\delta(Q,P)$ is a linear combination of morphisms of the form $U\times_X\delta$, where $U$ is a $(Q\times P\op,X)$-biset, i.e. a $(Q,P\times X)$-biset. So proving that $(ii)$ holds is equivalent to proving that $F(U\times_X\delta)=0$ for any $P$ and $Q$ and any such biset $U$. This will be a consequence of the next two lemmas, the first of which uses the following notation~: let $\tilde{U}$ be the $(X\times Q,P)$-biset equal to $U$ as a set, with double action defined by
$$\underbrace{(x,h)\mdot u\mdot g}_{{\rm in}~\tilde{U}}=\underbrace{hu(x^{-1},g)}_{{\rm in}~U}$$
for $x\in X$, $h\in Q$, $u\in U$ and $g\in P$. Then~:
\begin{enonce}{Lemma} \label{double}With this notation, let $T\in B(X\op)\cong B(\un\times X\op)$. So $\pi_Q(T)\in B\big(Q\times(Q\times X)\op\big)$, and $T\op\in B(X)$. Then 
$$U\times_X T\op=\pi_Q(T)\times_{(X\times Q)}\tilde{U}$$
in $B(Q\times P\op)$.
\end{enonce}
\pf I can suppose that $T$ is a right $X$-set. Then it is easily checked that the maps
\begin{eqnarray*}
(u,t)\in U\times_XT\op&\mapsto& \big((t,1),u\big)\in \pi_Q(T)\times_{(X\times Q)}\tilde{U}\\
\big((t,q),u\big)\in \pi_Q(T)\times_{(X\times Q)}\tilde{U}&\mapsto &(qu,t)\in U\times_XT\op
\end{eqnarray*} 
are well defined mutual inverse isomorphisms of $(Q,P)$-bisets.\findemo
\begin{enonce}{Lemma}\label{delta nul} If $F$ is a rational $p$-biset functor, then the map 
$$F(\delta\op)~:F(X)\to F(\un)$$
is equal to 0.
\end{enonce}
\pf Recall that $\delta=(X/I-X/IZ)-(X/J-X/JZ)$, i.e.
$$\delta\op=(I\dom X-IZ\dom X)-(J\dom X-JZ\dom X)\mvirg$$
where $Z$ is the center of $X$, and where $I$ and $J$ are non central subgroups of order $p$ in $X$, not conjugate in $X$.\par
Let $\mathcal{G}$ the set of subgroups of $X$ defined by
$$\mathcal{G}=\{Y\subseteq X\mid |X:Y|\leq p\}\sqcup \{I\}\mpoint$$
Then $\mathcal{G}$ is a genetic basis of $X$~: indeed the set $\{Y\subseteq X\mid |X:Y|\leq p\}$ is obtained by inflation from the unique genetic basis of the elementary abelian group $X/Z$. This set contains $p+2$ subgroups. Moreover, there is a unique faithful irreducible rational representation of $X$ (because there are $p+3$ conjugacy classes of cyclic subgroups in $X$), and it is easy to check that $I$ is a genetic subgroup of $X$ corresponding to this representation~: clearly $I\cap Z=\un$, and $N_X(I)/I=IZ/I$ is cyclic of order $p$. Moreover $\hat{I}=IZ$, and if $x\in X$ is such that $I^x\cap IZ\subseteq I$, then $I^x\subseteq I$ because $IZ$ is a normal subgroup of $X$. So $I^x=I$, and $I$ is a genetic subgroup of $X$, and $\mathcal{G}$ is a genetic basis of $X$.\par
Let $Y\in\mathcal{G}$. If $Y=P$, then $\gamma_Y$ is an $(X,X)$-biset with one element. So $T\dom X\times\gamma_P$ is a set of cardinality~1, for any subgroup $T$ of $X$. Hence in $B(X\op)$
$$\delta\op\times_X\gamma_P=0\mpoint$$
Now if $Y$ has index $p$ in $X$, then $\gamma_Y=X/Y-X/X$, thus
$$\delta\op\times_X\gamma_Y=\delta\op\times_X X/Y\mpoint$$
But $T\dom X\times_XX/Y\cong TY\dom X$, so $T\dom X\times_XX/Y\cong TZ\dom X\times_XX/Y$, for any subgroup $T$ of $X$, since $Y\normal X$ and $Y\supseteq Z$. Thus again
$$\delta\op\times_X\gamma_Y=O$$
in $B(X\op)$, in this case .\par
Finally if $Y=I$, then
$$\gamma_Y=(X/I\times_{IZ}I\dom X)-(X/IZ\times_{IZ}IZ\dom X)=(X/I\times_{IZ}I\dom X)-X/IZ$$
since $IZ\normal X$. By the previous argument, it follows that
$$\delta\op\times_X\gamma_Y=\delta\op\times_X(X/I\times_{IZ}I\dom X)\mpoint$$
Now $X/I\times_{IZ}I\dom X\cong (X\times X)/\Delta_{IZ,Z}$, and by the Mackey formula for the product of bisets (Proposition~1 of~\cite{doublact}), for any subgroup $T$ of $X$
$$T\dom X\times_X(X\times X)/\Delta_{IZ,Z}=\sum_{x\in T\dom X/IZ}T_x\dom X\mvirg$$
where 
\begin{eqnarray*}
T_x&=&\{a\in X\mid \exists t\in T,\;(t^x,a)\in\Delta_{IZ,I}\}\\
&=&\{a\in X\mid \exists t\in T,\;t^x\in IZ,\;t^xa^{-1}\in I\}\\
&=&I(T^x\cap IZ)\mpoint
\end{eqnarray*}
Moreover $T\dom X/IZ=TIZ\dom X$ since $IZ\normal X$. Hence~:
\begin{eqnarray*}
I\dom X\times_X (X\times X)/\Delta_{IZ,Z}&=&I\dom X+IZ\dom X\\
IZ\dom X\times_X (X\times X)/\Delta_{IZ,Z}&=&2\,IZ\dom X\\
J\dom X\times_X (X\times X)/\Delta_{IZ,Z}&=&I\dom X\\
JZ\dom X\times_X (X\times X)/\Delta_{IZ,Z}&=&IZ\dom X\\
\end{eqnarray*}
Thus
$$\delta\op\times_X\gamma_Y=(I\dom X+IZ\dom X)-2\,IZ\dom X-I\dom X+IZ\dom X=0$$
again, in this case.\par
But since $F$ is rational, the sum $\sum_{Y\in\mathcal{G}}\limits F(\gamma_Y)$ is equal to the identity map of $F(X)$. Hence
$$F(\delta\op)=F(\delta\op)\circ \sum_{Y\in\mathcal{G}}F(\gamma_Y)=\sum_{Y\in\mathcal{G}}F(\delta\op\times_X\gamma_Y)=0$$
as was to be shown.  \findemo
\noindent{\bf End of the proof of Theorem~\ref{le mur}~:} it follows from Lemma~\ref{double} that
$$F(U\times_X\delta)=F\big(\pi_Q(\delta\op)\big)\circ F(\tilde{U})\mpoint$$
Now $F\big(\pi_Q(\delta\op)\big)=F_Q(\delta\op)$, and $F_Q$ is rational by Theorem~\ref{Yoneda-Dress rationnel}. Thus $F_Q(\delta\op)=0$ by Lemma~\ref{double}, and $F(U\times_X\delta)=0$, completing the proof of Theorem~\ref{le mur}.\findemo
\begin{rem}{Remark} \label{barkerrat} Consider the quotient category $\mathcal{C}_p^{lin}$ of $\mathcal{C}_p$, whose objects are finite $p$-groups, and morphisms are defined by
$$\Hom_{\mathcal{C}_p^{lin}}(P,Q)=R_\Q(Q\times P\op)\mpoint$$
Then a rational biset functor for L. Barker's definition is just an additive functor from $\mathcal{C}_p^{lin}$ to $\mathcal{A}b$.\par
When $p$ is odd, it follows from Corollary~\ref{th sur K} that the category $\mathcal{C}_p/\delta$ and~$\mathcal{C}_p^{lin}$ coincide.
So in this case, this definition of a rational biset functor is equivalent to Definition~\ref{defrat}. However if $p=2$, those two definitions are not equivalent, as will be shown in Section~\ref{conoyau exp}.
\end{rem}
\begin{rem}{Remark} If $F$ is rational, then $F_Q$ is also rational for any $Q$, thus $F_Q(\delta\op)=0$. Conversely, if $F_Q(\delta\op)=0$ for any $Q$, then the previous proof shows that $F$ is rational. This leads to the following~:
\end{rem}
\begin{enonce}{Corollary} \label{geometric}Let $S$ denote a Sylow $p$-subgroup of $PGL(3,\F_p)$. Then~$S$ acts on the set $\mathbb{P}$ of points and on the set $\mathbb{L}$ of lines of the projective plane over $\F_p$. \par
Let $F$ be a $p$-biset functor. Then $F$ is rational if and only if for any $p$-group $Q$, 
$$F(\mathbb{P}\op\times Q)=F(\mathbb{L}\op\times Q)$$
as maps from $F(S\times Q)$ to $F(Q)$, images by $F$ of the $(Q,S\times Q)$-bisets $\mathbb{P}\op\times Q$ and $\mathbb{L}\op\times Q$ respectively.
\end{enonce}
\pf Indeed $S$ is isomorphic to the group $X$, so I can suppose $S=X$. Now it is easy to check that $S$ has exactly one fixed point on $\mathbb{P}$ and on $\mathbb{L}$, and that for suitable choices of $I$ and $J$, the permutation representations of~$S$ on $\mathbb{P}$ and $\mathbb{L}$ are respectively
$$\mathbb{P}\cong S/S\sqcup S/IZ\sqcup S/J\ressort{1cm}\mathbb{L}\cong S/S\sqcup S/JZ\sqcup S/I\mvirg$$
so that $\mathbb{P}-\mathbb{L}=\delta$ in $B(S)$.\findemo
\begin{rem}{Remark} Let $A$ be an injective $\Z$-module. If $F$ is a rational $p$-biset functor, then $F^0=\Hom(F,A)$ is rational by Proposition~7.4 of~\cite{bisetsections}. Conversely, if $F^0$ is rational, then $(F^0)^0$ is rational, and since the canonical map $F\to (F^0)^0$ is injective (because $A$ is an injective $\Z$-module), it follows that $F$ is rational. This duality argument shows that a $p$-biset functor $F$ is rational if and only if for any $p$-group $Q$, one has that
$$F(\mathbb{P}\times Q)=F(\mathbb{L}\times Q)\mvirg$$
as maps from $F(Q)$ to $F(S\times Q)$, images by $F$ of the $(S\times Q,Q)$-bisets $\mathbb{P}\times Q$ and $\mathbb{L}\times Q$ respectively. Equivalently $F_Q(\delta)=0$.
\end{rem}
\begin{enonce}{Proposition} Let $F$ be a $p$-biset functor.
\begin{enumerate}
\item The correspondence sending a $p$-group $P$ to the image of the map  $F_P(\delta\op)~:F(X\times P)\to F(P)$ defines a biset subfunctor of $F$, and the quotient functor $F^{rat}$ is the largest rational quotient of $F$. Thus
$$F^{rat}(P)=F(P)/\Im\, F_P(\delta\op)\mpoint$$
\item The correspondence sending a $p$-group $P$ to the kernel of the map $F_P(\delta): F(P)\to F(X\times P)$ defines a biset subfunctor $F_{rat}$ of $F$, which is the largest rational subfunctor of $F$. Thus
$$F_{rat}(P)=\Ker\,F_P(\delta)\mpoint$$
\end{enumerate}
\end{enonce}
\pf Let $P$ and $Q$ be finite $p$-groups, and let $U$ be a finite $(Q,P)$-biset. It is easy to check that the following diagram 
$$\carre{F(X\times P)}{F_P(\delta\op)}{F(P)}{F(X\times U)}{F(U)}{F(X\times Q)}{F_Q(\delta\op)}{F(Q)}$$
is commutative. It shows that the correspondence $P\mapsto \Im\, F_P(\delta\op)$ defines a subfunctor of $F$. Now if $N$ is a subfunctor of $F$, then $F/N$ is rational if and only if for any $p$-group $P$, the map $(F/N)_P(\delta\op)$ is equal to 0. This is equivalent to $\Im\, F_P(\delta\op)\subseteq N(P)$, so $F^{rat}$ is the largest rational quotient of~$F$.\par
Similarly, for Assertion 2, the diagram
$$\carre{F(P)}{F_P(\delta)}{F(X\times P)}{F(U)}{F(X\times U)}{F(Q)}{F_Q(\delta)}{F(X\times Q)}$$
is commutative, showing that $F_{rat}$ is a subfunctor of $F$. And if $M$ is a subfunctor of $F$, then $M$ is rational if and only if for any $p$-group~$P$, the map~$M_P(\delta)$ is equal to 0, i.e. if $M(P)\subseteq F_{rat}(P)$.\findemo
\begin{enonce}{Corollary} The functor $B/B_\delta$ is the largest rational quotient of $B$. If $p$ is odd, it is equal to $R_\Q$.
\end{enonce}
\pf Indeed $B/B_\delta$ is rational by Theorem~\ref{brat}. Conversely, the image of $B_\un(\delta): B(\un)\to B(X)$ contains $\delta$, so if $B'$ is a subfunctor of $B$ such that $B/B'$ is rational, then $B'\supseteq B_\delta$. Finally $B/B_\delta=R_\Q$ if $p$ is odd, by Corollary~\ref{th sur K}.\findemo
\section{The functor $K/B_\delta$\label{conoyau exp}}
Since $B/K\cong R_\Q$, there is an exact sequence of $p$-biset functors
$$0\to K/B_\delta\to B/B_\delta\to R_\Q\to 0\mpoint$$
Moreover, the functor $B/B_\delta$ is rational by Theorem~\ref{brat}. It follows that its subfunctor $K/B_\delta$ is also rational, and Corollary~\ref{th sur K} shows that this functor is equal to zero if $p$ is odd. This section describes the structure of $K/B_\delta$ for $p=2$, but the methods and results exposed here also hold for $p$ odd, and for this reason $p$ still denotes an arbitrary prime number here.\par
Since $K/B_\delta$ is rational, then for any $p$-group $P$ and any genetic basis $\mathcal{G}$ of~$P$, the map
$$\dirsum_{Q\in\mathcal{G}}\limits\Indinf_{N_P(Q)/Q}^P~:~\dirsum_{Q\in\mathcal{G}}\limits\partial (K/B_\delta)\big(N_P(Q)/Q\big)\to (K/B_\delta)(P)$$
is an isomorphism. Hence the evaluation of $K/B_\delta$ at $P$ is known if the groups $\partial (K/B_\delta)\big(N_P(Q)/Q\big)$ are known, when $Q$ is a genetic subgroup of $P$. Since $N_P(Q)/Q$ has normal $p$-rank~1 in this case, this relies on the following lemma~:
\begin{enonce}{Lemma} Let $R$ be a $p$-group of normal $p$-rank~1. If $R$ is non trivial, denote by $Z$ the unique subgroup of order $p$ in the center of $R$. Then~:
\begin{enumerate}
\item If $R$ is cyclic or generalized quaternion, then $\partial K(R)=\zero$, hence $\partial (K/B_\delta)(R)=\zero$.
\item If $R$ is semidihedral, then $\partial K(R)$ is free of rank one, generated by the element
$$ -2(R/W-R/WZ)+(R/\un -R/Z)=\Ind_{WZ}^R\varepsilon_{WZ}\mvirg$$
where $W$ is a non central subgroup of order~2 in $R$. In particular $\partial (K/B_\delta)(R)=\zero$.
\item If $R$ is dihedral, then $\partial K(R)$ is free of rank two, generated by the elements
$$ -2(R/W-R/WZ)+R/\un -R/Z=\Ind_{WZ}^R\varepsilon_{WZ}\mvirg\;\hbox{and}$$
$$\delta_R=(R/W-R/WZ)-(R/W'-R/W'Z)\mvirg$$
where $W$ and $W'$ are non central subgroups of order~2 in $R$, which are not conjugate in $R$. In this case $\partial (K/B_\delta)(R)\cong \Z/2\Z$, generated by the image of $\delta_R$.
\end{enumerate}
\end{enonce}
\pf If $R$ is trivial, then $K(R)=\zero$, so the result holds. And if $R$ is non trivial, then the idempotent $f_\un^R$ is equal to $R-R/Z$, so for any subgroup $S$ of~$R$
$$f_\un^RR/S=R/S-R/SZ\mvirg$$
and this is zero if $S\supseteq Z$. \par
Now if $R$ is cyclic or generalized quaternion, and if $S\not\supseteq Z$, then $S=\un$. Thus $\partial B(R)$ is free of rank one in this case, generated by $R/\un -R/Z$. Since $|R/\un-R/Z|=|R|/p$, and since $|Y|=0$ for any $Y\in K(R)$, it follows that $\partial K(R)=\zero$ in this case, proving Assertion~1.\par
If $R$ is semidihedral, and $S\not\supseteq Z$, then either $S=\un$ or $S=W$, up to conjugation in $R$. So $\partial B(R)$ is free of rank~2, generated by $R/\un-R/Z$ and $R/W-R/WZ$. Since $|R/\un-R/Z|=|R|/2$ and $|R/W-R/WZ|=|R|/4$, any element in $\partial K(R)$ is a multiple of
$$u=-2(R/W-R/WZ)+(R/\un-R/Z)\mvirg$$
which is equal to $\Ind_{WZ}^R\varepsilon_{WZ}$. Since $\varepsilon_{WZ}\in K(WZ)$, it follows that $\partial K(R)$ is free of rank one, generated by $u$. Assertion~2 follows, since moreover $u\in K_\varepsilon(R)$ and $K_\varepsilon(R)\subseteq K_\delta(R)$.\par
Finally if $R$ is dihedral, and $S\not\supseteq Z$, then either $S=\un$, or $S$ has order~2, and $S\neq Z$. So $\partial B(R)$ is free of rank~3, generated by the elements $R/\un-R/Z$, $R/W-R/WZ$, and $R/W'-R/W'Z$, where $W$ and $W'$ are representatives of the two conjugacy classes of non central subgroups of order~2 in~$R$. Looking at cardinalities shows that $\partial K(R)$ is contained in the submodule of $\partial B(R)$ generated by the elements
$$u=-2(R/W-R/WZ)+(R/\un-R/Z)\mvirg\; \hbox{and}$$
$$\delta_R=(R/W-R/WZ)-(R/W'-R/W'Z)\mpoint$$
Conversely, since $u=\Ind_{WZ}^R\varepsilon_{WZ}$, it follows that $u\in K(R)$. But $\delta_R\in K(R)$ also, by Remark~6.10 of~\cite{dadegroup}, or by direct computation of the character of~$\delta_R$. Hence $\partial K(R)$ is equal to the module generated by $u$ and $\delta_R$. Since moreover $u\in K_\varepsilon(R)\subseteq K_\delta(R)$ as before, and since
$$2\delta_R=\Ind_{W'Z}^R\varepsilon_{W'Z}-\Ind_{WZ}^R\varepsilon_{WZ}\in K_\varepsilon(R)\subseteq K_\delta(R)\mvirg$$ 
it follows that $\partial K/B_\delta(R)$ is either trivial or of order~2, generated by the image of $\delta_R$. The next lemma will show that $\delta_R\notin B_\delta(R)$, completing the proof of Assertion~3.\findemo
\begin{enonce}{Lemma} Suppose $p=2$, and let $R$ be a dihedral group of order at least 16. Then
$\delta_R\notin B_\delta(R)$.
\end{enonce}
\pf Suppose that $\delta_R\in B_\delta(R)$. Then there exists $\varphi\in \Hom_{{\mathcal{C}_p}}(D_8,R)$ such that $\delta_R=\varphi(\delta)$. In particular there exists some $(R,D_8)$-biset $U$ such that the coefficient of $R/W$ in the expansion of $U\times_{D_8}\delta$ in the canonical basis of $B(R)$ is odd, and I can assume that $U$ is a transitive biset, of the form $U=(R\times D_8)/L$, for some subgroup $L$ of $R\times D_8$. Such a transitive biset factors as
$$\Indinf_{T/S}^R\Iso_{B/A}^{T/S}\Defres_{B/A}^{D_8}\mvirg$$
where $(T,S)$ is a section of $R$ and $(B,A)$ is a section of $D_8$ such that the factor groups $T/S$ and $B/A$ are isomorphic. Now since $\delta\in K(D_8)$, and since $K(C)=\zero$ for a cyclic group $C$, it follows that $\Defres_{B/A}^{D_8}\delta=0$, except if $(B,A)=(D_8,\un)$ or if $B/A$ is elementary abelian of rank~2. Moreover any proper deflation of $\delta$ is equal to~0, hence in the latter case $B$ cannot be equal to $D_8$. So $B$ is one of the two elementary abelian subgroups of index 2 in $D_8$, and $A=\un$. In this case $\Defres_{B/A}^{D_8}\delta$ is equal to the restriction of $\delta$ to $B$, i.e. to $\pm\varepsilon_B$.\par
Thus there is a section $(T,S)$ of $R$, such that $T/S$ is isomorphic to $D_8$ or $(C_2)^2$, such that $U\times_{D_8}\delta=\Indinf_{T/S}^Ru$, where $u=\delta_{T/S}$ if $T/S\cong D_8$, and $u=\varepsilon_{T/S}$ if $T/S\cong (C_2)^2$. The coefficient of $R/W$ in $\Indinf_{T/S}^Ru$ is equal to 0, except if $S\subseteq W^r\subseteq T$ for some $r\in R$. Thus either $S=\un$ or $S$ is conjugate to $W$. But $N_R(W)/W$ has order 2, so it cannot contain a group isomorphic to $T/S$. It follows that $S=\un$. In other words $U\times_{D_8}\delta=\Ind_T^Ru$, where $T$ is a subgroup of $R$ which is isomorphic to $D_8$, and in this case $u=\delta_T$, or isomorphic to $(C_2)^2$, and $u=\varepsilon_T$. In the first case
$$\Ind_T^R\delta=(R/I-R/IZ)-(R/J-R/JZ)\mvirg$$
where $I$ and $J$ are non central subgroups of order 2 of $D_8$. Since $|R|\geq 16$, the groups $I$ and $J$ are conjugate in $R$, so $\Ind_T^R\delta=0$. In the second case $T\supseteq Z$, and
$$\Ind_T^R\varepsilon_T=R/\un -(R/A+R/B+R/Z)+2R/T\mvirg$$
where $A$, $B$ and $Z$ are the subgroups of order 2 of $T$. But $A$ and $B$ are conjugate in $R$, so in fact
$$\Ind_T^R\varepsilon_T=R/\un -(2R/A+R/Z)+2R/T\mvirg$$
The coefficient of $R/W$ in this expression is equal to -2 or 0, according to $A$ being conjugate to $W$ in $R$ or not. In any case, it is even, and this completes the proof.\findemo
\begin{enonce}{Corollary} \label{base K mod delta}Let $P$ be a $p$-group, and let $\mathcal{G}$ be a genetic basis of~$P$. Then $K(P)/B_\delta(P)$ is a vector space over $\F_2$, with basis the images of the elements $\Indinf_{N_P(Q)/Q}^P\delta_{N_P(Q)/Q}$, where $Q$ runs through the elements of $\mathcal{G}$ for which $N_P(Q)/Q$ is dihedral.
\end{enonce}
So $K(P)/B_\delta(P)\cong (\Z/2\Z)^d$, where $d$ is the number of isomorphism classes of rational irreducible representations of $P$ of dihedral type. This result suggest a little more, in view of the following result~:
\begin{enonce}{Theorem} If $P$ is a $p$-group, denote by $B^\times(P)$ the group of units of its Burnside ring $B(P)$.
\begin{enumerate}
\item {\rm[\cite{burnsideunits} Theorem~7.2 and Corollary~7.3]} The correspondence $P\mapsto B^\times(P)$ is a $p$-biset functor, which is isomorphic to a subfunctor of $\F_2R_\Q^*$. In particular, it is a rational $p$-biset functor. 
\item {\rm[\cite{burnsideunits} Theorem~8.4]} The functor $B^\times$ is uniserial, and the poset of its proper subfunctors
$$\zero\subset L_0\subset L_1\ldots \subset L_n\subset \ldots$$
is such that $L_0$ is isomorphic to the simple biset functor $S_{\un,\F_2}$, and $L_i/L_{i-1}$, for $i\geq 1$, is isomorphic to the simple functor $S_{D_{2^{i+3}},\F_2}$. 
\item {\rm[\cite{burnsideunits} Theorem~8.4]} There is a map of $p$-biset functors ${\rm exp}: B\to B^\times$, whose image is equal to the socle $L_0$ of $B^\times$. So for any $p$-group $P$, the $\F_2$-dimension of the cokernel of ${\rm exp}_P: B(P)\to B^\times(P)$ is equal to the number of isomorphism classes of rational irreducible representations of $P$ of dihedral type.
\end{enumerate}
\end{enonce}
This suggest that the biset functors $K/B_\delta$ and $\Coker\;{\rm exp}=B^\times/{\rm soc}(B^\times)$ might be isomorphic. And this is indeed the case~:
\begin{enonce}{Theorem} \label{coker}The functor $K/B_\delta$ is isomorphic to the cokernel of the map ${\rm exp}~: B\to B^\times$. So $K=B_\delta$ if $p\neq 2$, and if $p=2$, the functor $K/B_\delta$ is uniserial. Its lattice of proper subfunctors is
$$\zero=\sur{K}_0\subset \sur{K}_1\subset \sur{K}_2\subset\ldots\subset \sur{K}_n\subset\ldots$$
where $\sur{K}_n$, for $n\geq 1$, is generated by the image of $\delta_{D_{2^{n+3}}}\in K(D_{2^{n+3}})$ in $(K/B_\delta)(D_{2^{n+3}})$. Moreover $\sur{K}_n/\sur{K}_{n-1}\cong S_{D_{2^{n+3}},\F_2}$, for $n\geq 1$. 
\end{enonce}
\pf First recall that if $R$ is a $p$-group of normal $p$-rank~1, then $R$ has a unique faithful rational irreducible representation $\Phi_R$, up to isomorphism, by Proposition~3.7 of~\cite{fonctrq}. For any $p$-group $P$, the group $R_\Q(P)$ is a free abelian group with basis the set $\Irr_\Q(P)$ of isomorphism classes of irreducible rational representations of $P$. So the dual group $R_\Q^*(P)$ has a dual basis $(V^*)_{V\in \Irr_\Q(P)}$, where the value of $V^*$ on some $\Q P$-module $M$ is equal to the multiplicity of~$V$ as a direct summand of $M$ (with the slight abuse of notation identifying an irreducible representation $V$ with its isomorphism class). Denote by $\sur{V}^*$ the image of $V^*$ in $\F_2R_\Q^*$.\par
Also set $K_\flat=K/B_\delta$, and $B^\times_\flat=B^\times/L_0$.
With this notation, and identifying moreover the functor $B^\times$ with its isomorphic image inside $\F_2R_\Q^*$, the subfunctor $L_n$ of $B^\times$ is generated by the element $\sur{\Phi}^*_{D_{2^{n+3}}}\in \F_2R_\Q^*(D_{2^{n+3}})$, for $n\geq 1$, and the functor $L_0$ is generated by the element $\sur{\Phi}_\un^*\in\F_2R_\Q^*(\un)$. So if $P$ is a $p$-group, and if $\mathcal{G}$ is a genetic basis of $P$, then the group $B^\times_\flat$ has a basis over $\F_2$ consisting of the elements $\sur{V}_S^*=\Indinf_{N_P(S)/S}^P\sur{\Phi}_{N_P(S)/S}^*$, where $S$ runs through the set $\mathcal{D}$ of elements of~$\mathcal{G}$ for which $N_P(S)/S$ is dihedral.\par
Moreover, if $R$ is a $p$-group of normal $p$-rank~1, then $\partial B^\times(R)$ is equal to~$\zero$, unless $R$ is dihedral or cyclic of order at most~2, and in each of these cases $\partial B^\times(R)\cong\F_2$ (see Yal\c c\i n \cite{yalcin} Lemma~4.6 and Lemma~5.2, or Corollary~5.9 of~\cite{burnsideunits}). But if $|R|\leq 2$, then $S_{D_{2^{n+3}},\F_2}(R)=\zero$ for $n\geq 1$, thus $B^\times(R)=L_0(R)$, and it follows that $\partial B^\times_\flat(R)=\zero$, except if $R$ is dihedral. In this case since $\partial S_{R,\F_2}(R)\cong\F_2\neq\zero$, it follows that $\partial B^\times_\flat(R)\cong \F_2$, generated by the image of the element $\sur{\Phi}_R^*$. It will be useful for further computation to notice that $\sur{\Phi}_R^*(\Q R/T)$ is equal to 2 if $T=\un$, to 1 if $T$ is a non central subgroup of order 2 of~$R$, and to 0 otherwise.\par
On the other hand, by Corollary~\ref{base K mod delta}, the group $K_\flat(P)$ has an $\F_2$-basis consisting of the images $\sur{d}_S$ in $K_\flat(P)$ of the elements $\Indinf_{N_P(S)/S}^P\delta_{N_P(S)/S}$, for $S\in\mathcal{D}$. So there is a unique linear map
$$\psi_{P,\mathcal{G}}: K_\flat(P)\to B^\times_\flat(P)$$
sending $\sur{d}_S$ to $\sur{V}_S^*$, for $S\in\mathcal{D}$, and this map is an isomorphism. In other words, the map $\psi_{P,\mathcal{G}}$ is equal to the composition of isomorphisms
\begin{equation}\label{decomposition}K_\flat(P)\stackrel{\beta}{\to}\dirsum_{S\in\mathcal{G}}\partial K_\flat\big(N_P(S)/S\big)\stackrel{\gamma}{\to}\dirsum_{S\in\mathcal{G}}\partial B^\times_\flat\big(N_P(S)/S\big)\stackrel{\alpha}{\to} B^\times_\flat(P)\mvirg
\end{equation}
where $\alpha=\dirsum_SB^\times_\flat(a_S)$ and $\beta=\dirsum_SK_\flat(b_S)$, and where $\gamma$ is the direct sum of the unique isomorphisms $\gamma_{N_P(S)/S}~:\partial K_\flat\big(N_P(S)/S\big)\to \partial B^\times_\flat\big(N_P(S)/S\big)$~: this makes sense because the groups $\partial K_\flat\big(N_P(S)/S\big)$ and $B^\times_\flat\big(N_P(S)/S\big)$ are both cyclic of order 2, when $S\in\mathcal{D}$, or both trivial if $S\in\mathcal{G}-\mathcal{D}$. 
\begin{enonce}{Lemma}\label{fonctoriel} Let $P$ and $Q$ be $p$-groups. Let $\mathcal{G}$ be a genetic basis of $P$, and $\mathcal{H}$ be a genetic basis of $Q$. Then if $\varphi\in\Hom_{\mathcal{C}_p}(P,Q)$, the diagram
$$\carre{K_\flat(P)}{\psi_{P,\mathcal{G}}}{B^\times_\flat(P)}
{K_\flat(\varphi)}{B^\times_\flat(\varphi)}
{K_\flat(Q)}{\psi_{Q,\mathcal{H}}}{B^\times_\flat(Q)}$$
is commutative.
\end{enonce}
\pf Using decomposition~\ref{decomposition} both for $P$ and $Q$, this diagram gives
$$\begin{array}{ccccccc}
K_\flat(P)&\!\!\!\stackrel{\beta}{\to}&\!\!\!\dirsum_{S\in\mathcal{G}}\limits\partial K_\flat\big(N_P(S)/S\big)&\!\!\!\stackrel{\gamma}{\to}&\!\!\!\dirsum_{S\in\mathcal{G}}\limits\partial B^\times_\flat\big(N_P(S)/S\big)&\!\!\!\stackrel{\alpha}{\to}&\!\!\! B^\times_\flat(P)\\
\flv{K_\flat(\varphi)}{}&&&&&&\flv{}{B^\times_\flat(\varphi)}\\
K_\flat(Q)&\!\!\!\stackrel{\beta'}{\to}&\!\!\!\dirsum_{T\in\mathcal{H}}\limits\partial K_\flat\big(N_Q(T)/T\big)&\!\!\!\stackrel{\gamma}{\to}&\!\!\!\dirsum_{T\in\mathcal{H}}\limits\partial B^\times_\flat\big(N_Q(T)/T\big)&\!\!\!\stackrel{\alpha'}{\to}&\!\!\! B^\times_\flat(Q)\\
\end{array}
$$
where the maps $\alpha'$ and $\beta'$ are the analogues of $\alpha$ and $\beta$  for the group~$Q$ and its genetic basis $\mathcal{H}$. Showing that this diagram is commutative amounts to showing that
$$B^\times_\flat(\varphi)\circ \alpha\circ\gamma\circ\beta=\alpha'\circ\gamma'\circ\beta'\circ K_\flat(\varphi)\mvirg$$
or equivalently, that
$$\alpha'^{-1}\circ B^\times_\flat(\varphi)\circ \alpha\circ\gamma=\gamma\circ\beta'\circ K_\flat(\varphi)\circ\beta^{-1}\mpoint$$
Now $\beta^{-1}=\dirsum_{S\in\mathcal{G}}\limits K_\flat(a_S)$ and $\alpha'^{-1}=\dirsum_{T\in\mathcal{H}}\limits B^\times_\flat(b_T)$. So this equality is equivalent~to
\begin{equation}\label{à démontrer}
B^\times_\flat(b_T\varphi a_S)\circ\gamma_{N_P(S)/S}=\gamma_{N_Q(T)/T}\circ K_\flat(b_T\varphi a_S)\mvirg
\end{equation}
for any $S\in\mathcal{G}$ and any $T\in\mathcal{H}$. Now 
$$b_T\varphi a_S\in f_\un^{N_Q(T)/T}\Hom_{\mathcal{C}_p}\big(N_P(S)/S,N_Q(T)/T\big)f_\un^{N_P(S)/S}\mvirg$$
so equality~(\ref{à démontrer}), hence Lemma~\ref{fonctoriel}, follow from the next lemma.\findemo
\begin{enonce}{Lemma}\label{cas prn1} Let $M$ and $N$ be $p$-groups of normal $p$-rank~1, and let $f\in f_\un^N\Hom_{\mathcal{C}_p}(M,N)f_\un^M$, then
$$B_\flat^\times(f)\circ\gamma_M=\gamma_N\circ K_\flat(f)\mpoint$$
\end{enonce}
\pf Both sides are maps from $\partial K_\flat(M)$ to $\partial B_\flat^\times(N)$, so if one of the groups~$M$ or~$N$ is not dihedral, both sides are equal to zero. So I can suppose that $M$ and $N$ are dihedral. I can also suppose that $f=f_\un^NUf_\un^P$, for some transitive $(N,M)$-biset $U$ of the form $(N\times M)/L$, where $L$ is a subgroup of $N\times M$. If $f=0$, there is nothing to prove, and otherwise $k_1(L)\cap Z(N)=\un$ and $k_2(L)\cap Z(M)=\un$. It follows that $k_1(L)$ and $k_2(L)$ are trivial or non central of order~2. If one of these groups, say $k_1(L)$, is non central of order 2, then its normalizer has order 4, and contains $p_1(L)$. Hence the group $C=q(L)\cong p_1(L)/k_1(L)$ has order 1 or 2. But $f$ factors through~$C$ and both $K_\flat(C)$ and $B_\flat^\times(C)$ are equal to $\zero$. So again there is nothing to prove in this case.\par
Hence I can assume that $k_1(L)=k_2(L)=\un$. In this case, the morphism $(N\times M)/L$ is the composition of the restriction from $M$ to $p_2(L)$, followed by a group isomorphism from $p_2(L)$ to $p_1(L)$, followed by induction from $p_1(L)$ to $N$. In particular, it factors through $H=p_1(L)$. If this group is cyclic, or contained in a dihedral subgroup of order~8 of $M$, then $K_\flat(H)=B^\times_\flat(H)=\zero$ because $H$, has no genetic subgroup of dihedral type. So I can suppose that~$H$ is dihedral of order at least 16.\par
Now $\partial K_\flat(M)$ is generated by the image of the element
$$\delta_M=(M/W-M/WZ)-(M/W'-M/W'Z)\mvirg$$
where $Z$ is the center of $M$, and $W$ and $W'$ are non central subgroups of order~$2$, non conjugate in $M$. If $M'$ is a dihedral group of index 2 in $M$, containing $W$, then one checks easily that in $B(M')$
$$\Res_{M'}^M\delta_M=\delta_{M'}-\Ind_{W^xZ}^{M'}\varepsilon_{W^xZ}\mpoint$$
where $x\in M-M'$, so $\Res_{M'}^M\delta_M=\delta_{M'}$ in $K_\flat(M')$ since  $$\Ind_{W^xZ}^{M'}\varepsilon_{W^xZ}\in B_\varepsilon(M')\subseteq B_\delta(M')\mpoint$$
By induction, it follows that $\Res_H^M\delta_M=\delta_H$ in $K_\flat(H)$. This element is sent to $\delta_{p_1(L)}$ by the isomorphism from $H=p_2(L)$ to $p_1(L)$, hence finally
$$K_\flat(f)(\delta_M)=f_\un^N\Ind_{p_1(L)}^N\delta_{p_1(L)}$$
in $K_\flat(N)$. Now there are two cases~: either $p_1(L)=N$ and $K_\flat(f)(\delta_M)=\delta_N$, or $p_1(L)$ is a proper dihedral subgroup of $N$. In this case, the non central subgroups of order~2 of $p_1(L)$ are conjugate in $N$, so
$\Ind_{p_1(L)}^N\delta_{p_1(L)}=0$ in $K\big(p_1(L)\big)$.\par
So either $p_1(L)\neq N$, and $\gamma_N\circ K_\flat(f)(\delta_M)=0$, or $p_1(L)=N$, and $\gamma_N\circ K_\flat(f)(\delta_M)$ is equal to the generator $\sur{\Phi}_N^*$ of $\partial B_\flat^\times(N)$. \par
On the other hand $\gamma_M(\delta_M)$ is equal to the image of $\sur{\Phi}_M^*$ in $B^\times_\flat(M)$. The restriction of $\sur{\Phi}_M^*$ to the dihedral subgroup $H=p_2(L)$ is easily seen to be equal to $\sur{\Phi}_H^*$, hence it is mapped to $\sur{\Phi}_{p_1(L)}^*$ by the isomorphism $p_2(L)\to p_1(L)$. Finally
$$B_\flat^\times\circ\gamma_M(\delta_M)=\Ind_{p_1(L)}^N\sur{\Phi}_{p_1(L)}^*\mpoint$$
If $p_1(L)=N$, this is equal to $\sur{\Phi}_N^*$. Otherwise it is induced from a proper dihedral subgroup. It is easily checked in this case that the values of the linear form $\Ind_{p_1(L)}^N\Phi_{p_1(L)}^*$ are all even, so $\Ind_{p_1(L)}^N\sur{\Phi}_{p_1(L)}^*=0$ in this case. This completes the proof of Lemma~\ref{cas prn1}.\findemo
\noindent{\bf End of the proof of Theorem~\ref{coker}~:} First apply Lemma~\ref{fonctoriel} in the case $Q=P$, and $\varphi=\Id_P$~: this shows that if $\mathcal{G}$ and $\mathcal{H}$ are two genetic bases of~$P$, then $\psi_{P,\mathcal{G}}=\psi_{P,\mathcal{H}}$. In other words, the map $\psi_{P,\mathcal{G}}$ does not depend on the choice of the genetic basis $\mathcal{G}$ of $P$, and it can be denoted by $\psi_P$.\par
Now Lemma~\ref{fonctoriel} shows that these maps $\psi_P$ define a morphism of $p$-biset functors from $K_\flat$ to $B^\times_\flat$, and this morphism is an isomorphism. The other assertions in Theorem~\ref{coker} follow from the structure of the functor $B^\times$ (see Theorem~8.4 of~\cite{burnsideunits}).\findemo 
\begin{rem}{Remark} It follows that for $p=2$, the class $\mathcal{R}$ of additive functors $\mathcal{C}_p^{lin}\to \mathcal{A}b$ defined in Remark~\ref{barkerrat} is not closed under extensions~: indeed, in the short exact sequence
$$0\to K/K_\delta\to B/K_{\delta}\to R_Q\to 0\mvirg$$
the functors $K/K_\delta$ and $R_\Q$ are in $\mathcal{R}$, but the functor $R_\Q=B/K$ is clearly the largest quotient of $B$ with this property. So $B/K_\delta$ is not in the class $\mathcal{R}$.\par
On the other hand, one can give the following alternative proof of the fact that the class of rational biset functors is closed under extensions~: in a short exact sequence 
$$0\to M\to L\to N\to 0\mvirg$$
where $M$ and $N$ are rational $p$-biset functors, i.e. additive functors $\mathcal{C}_p/\delta\to\nolinebreak \mathcal{A}b$, a straightforward argument shows that the middle term $L$ is a functor from the quotient category $\mathcal{C}_p/(\delta\times \delta)$ of $\mathcal{C}_p$, to $\mathcal{A}b$. And it turns out that the categories $\mathcal{C}_p/(\delta\times \delta)$ and $\mathcal{C}_p/\delta$ coincide~: one can indeed find a subgroup $Y$ of $X^3=X\times X^2$ such that
$$(X^3/Y)\times_{X^2}(\delta\times \delta)=\delta\mpoint$$
A tedious computation shows that using Notation~\ref{rappeldade}, one can take
$$Y=\{(xy^{-1}\varphi(x),x,y)\mid (x,y)\in X^2,\;xy^{-1}\in IZ\}\mvirg$$
where $\varphi:X\to J$ is any surjective homomorphism with kernel $IZ$.
\end{rem}

\nopagebreak
\leftline{Serge BOUC}
\leftline{CNRS - LAMFA - Universit\'e de Picardie-Jules Verne}
\leftline{33 rue St Leu - 80039 - Amiens Cedex 1 - France}
\leftline{\tt email : serge.bouc@u-picardie.fr}

\end{document}